\documentclass[12pt]{article}
\usepackage{amsmath}
\usepackage{amssymb}
\usepackage[dvips]{graphicx}
\usepackage{epsfig}
\def\eps{\varepsilon}
\font\tencmmib=cmmib10 \skewchar\tencmmib '60
\newfam\cmmibfam
\textfont\cmmibfam=\tencmmib

\font\tenmsb=msbm10

\def\Bbb#1{\hbox{\tenmsb#1}}

\def\lessim{\ \lower4pt\hbox{$
\buildrel{\displaystyle <}\over\sim$}\ }
\def\gessim{\ \lower4pt\hbox{$\buildrel{\displaystyle >}
\over\sim$}\ }

\def\eps{{\varepsilon}}

\newtheorem{proposition}{Proposition}

\newtheorem{theorem}{Theorem}

\makeatletter
\@addtoreset{equation}{section}

\makeatother

\font\tencmmib=cmmib10 \skewchar\tencmmib '60
\newfam\cmmibfam
\textfont\cmmibfam=\tencmmib

\font\tenmsb=msbm10

\def\Bbb#1{\hbox{\tenmsb#1}}

\def\lessim{\ \lower4pt\hbox{$
\buildrel{\displaystyle <}\over\sim$}\ }
\def\gessim{\ \lower4pt\hbox{$\buildrel{\displaystyle >}
\over\sim$}\ }

\def\eps{\varepsilon}

\def\go0{\to 0}

\def\leftitem#1{\item{\hbox to\parindent{\enspace#1\hfill}}}

\def\sg{\sigma}

\def\sg2{\sigma^2}

\def\__{_{\infty}}

\begin{document}

\title{\Large Rademacher processes and bounding the risk of function learning}
\author{
V. Koltchinskii\thanks{Department of Mathematics and
Statistics, University of New Mexico, Albuquerque, NM 87131, USA. E-mail:
vlad@math.unm.edu}
\ and D. Panchenko\thanks{Department of Mathematics and
Statistics, University of New Mexico, Albuquerque, NM 87131, USA. E-mail:
panchenk@math.unm.edu}
\thanks{The research of V. Koltchinskii is partially supported by NSA Grant
MDA904-99-1-0031. The research of D. Panchenko was 
partially supported by Boeing Computer Services Grant 3-48181}
}
\maketitle

 \pagestyle{myheadings}
 \markboth{Rademacher processes and function learning}
 {V. Koltchinskii and D. Panchenko}

\begin{abstract}
We construct data dependent upper bounds on the risk in function learning
problems. The bounds are based on the local norms of the Rademacher process 
indexed by 
the underlying function class and they do not require prior knowledge about the
distribution of training examples or any specific properties of the function 
class. Using Talagrand's type concentration inequalities for empirical and 
Rademacher 
processes, we show that the bounds hold with high probability that decreases 
exponentially
fast when the sample size grows. In typical situations that are frequently 
encountered in 
the theory of function learning, the bounds give nearly optimal rate of 
convergence 
of the risk to zero.   
\end{abstract}

\section{Local Rademacher norms and bounds on the risk: main results}

Let $(S, {\cal A})$ be a measurable space and let ${\cal F}$ be a class of ${\cal A}$-measurable
functions from $S$ into $[0,1].$ Denote ${\cal P}(S)$ the set of all probability measures on
$(S,{\cal A}).$ Let $f_0\in {\cal F}$ be an unknown \it target function. \rm  Given a probability
measure $P\in {\cal P}(S)$ (also unknown), let $(X_1,\dots ,X_n)$ be an i.i.d. sample in $(S,{\cal A})$
with common distribution $P$ (defined on a probability space $(\Omega, \Sigma, {\Bbb P})$).  
In computer learning theory, the problem of estimating $f_0,$ based on the 
labeled sample $(X_1,Y_1),\dots , (X_n,Y_n),$ where $Y_j:=f_0(X_j),\ j=1,\dots ,n,$ is referred to as 
\it function learning problem. \rm The so called \it concept learning \rm is a special case
of function learning. In this case, ${\cal F}:=\{I_C: C\in {\cal C}\},$ where ${\cal C}\subset {\cal A}$
is called a class of concepts (see Vapnik (1998), Vidyasagar (1996), Devroye, Gy\" orfi and
Lugosi (1996) for the account on statistical learning theory). 
The goal of function learning is to find an estimate 
$\hat f_n := \hat f_n ((X_1,Y_1),\dots ,(X_n,Y_n))$ of the unknown target function such that the
$L_1$-distance between $\hat f_n$ and $f_0$ becomes small with high probability as soon as the
sample size becomes large enough. The $L_1$-distance $P|\hat f_n -f_0|$ is often called \it
the risk \rm (also the generalization, or prediction error) of the estimate $\hat f_n.$ A class ${\cal F}$
is called \it probably approximately correctly (PAC) learnable \rm iff for all $\eps >0$
$$
\pi_n ({\cal F};\eps):=\sup_{P\in {\cal  P}(S)} \sup_{f_0\in {\cal  F}}
{\Bbb P}\bigl\{P|\hat f_n-f_0|\geq \eps\bigr\}\to 0\ {\rm as}\ n\to\infty.
$$
The bounds on the probability $\pi_n ({\cal F};\eps)$ are of importance in the theory.
Such bounds allow one to determine the quantity 
$$
N_{\cal  F}(\eps ;\delta):=\inf \{n: \pi_n ({\cal F};\eps)\leq \delta\},
$$
which is called \it the sample complexity of learning. \rm 
Unfortunately, a bound that is uniform in the class of all
distributions ${\cal P}(S)$ is not necessarily tight for a 
particular distribution $P$ and often such a bound does
not provide a reasonable estimate of the minimal sample
size needed to achieve certain accuracy of learning in
the case of a particular $P.$

A natural approach to the function learning problem (in the case
when $f_0\in {\cal F}$) is to find $\hat f_n\in {\cal F}$ such that 
$\hat f_n (X_j)=f_0(X_j)=Y_j$ for all $j=1,\dots ,n.$ 
In learning theory, such an estimate $\hat f_n$ is called \it consistent \rm (this notion should not be confused
with consistency in statistical sense).

We construct below a data dependent bound on the risk 
of a consistent estimate $\hat f_n.$ 
More precisely, given $\delta>0,$ we define a quantity
$$
\hat \beta_n ({\cal F};\delta) = \hat \beta_n ({\cal F};\delta; (X_1,Y_1),\dots, 
(X_n,Y_n))
$$
such that for any consistent estimate $\hat f_n$ 
\begin{equation}
\sup_{P\in {\cal  P}(S)}\sup_{f_0\in {\cal  F}}{\Bbb P}\bigl\{P|\hat f_n-f_0|\geq \hat \beta_n ({\cal F};\delta)\bigr\}
\leq \delta .
\end{equation}
We'll consider a couple of important examples in which the bound we suggest gives nearly optimal
rate of convergence of the risk to $0$ as the sample size tends to infinity.

Given a class ${\cal G}$ of ${\cal A}$-measurable
functions from $S$ into $[0,1]$ with $0\in {\cal G},$
let $\hat {\cal G}_n$ denote the restriction of the class ${\cal G}$ on the sample $(X_1,\dots ,X_n).$ 
Consider a quantity 
$$
\hat \gamma_n ({\cal G};\delta) = 
\hat \gamma_n (\hat {\cal G}_n;\delta; X_1,\dots, X_n)
$$
such that the bound 
$$
\sup_{P\in {\cal  P}(S)}{\Bbb P}\bigl\{P\hat g_n\geq \hat \gamma_n ({\cal G};\delta)\bigr\}
\leq \delta 
$$
holds for any class ${\cal G}$ and for any function 
$\hat g_n\in {\cal G}$ satisfying the conditions  
$\hat g_n(X_j)=0$ for all $j=1,\dots ,n.$

Define 
$${\cal F}(f_0):=\{|f-f_0|: f\in {\cal F}\}$$
(note that the values of the functions from
this class are known on the sample $(X_1,\dots ,X_n)$)  
and 
\begin{eqnarray*}
{\hat {\cal F}}_n(f_0)
&:=&
\{(|f-f_0|(X_j): 1\leq j\leq n): f\in {\cal F}\}
\\
&=&
\{(|f(X_j)-Y_j|: 1\leq j\leq n): f\in {\cal F}\}.
\end{eqnarray*} 
If $\hat f_n$ is a consitent estimate, then the function
$\hat g_n:=|\hat f_n - f_0|\in {\cal F}(f_0)$
satisfies the condition 
$\hat g_n(X_j)=0$ for all $j=1,\dots ,n.$ 
Then, clearly, for any consistent estimate $\hat f_n,$
$$
\sup_{P\in {\cal  P}(S)}\sup_{f_0\in {\cal  F}}{\Bbb P}\bigl\{P|\hat f_n-f_0|\geq \hat \gamma_n ({\cal F}(f_0);\delta)\bigr\}
\leq \delta .
$$
Therefore if one defines (for $Y_j=f_0(X_j)$) 
$$
\hat \beta_n ({\cal F};\delta; (X_1,Y_1),\dots, (X_n,Y_n)):=
\hat \gamma_n ({\hat {\cal F}}_n (f_0);\delta; X_1,\dots, X_n),
$$
then (1.1) holds. 

These considerations show that the problem can always be reduced to the case $f_0\equiv 0.$ 
To simplify the notations, we make this assumption in what follows. 

We also assume 
for simplicity that ${\cal F}$ is a countable class of functions. This condition
can be easily replaced by standard measurability assumptions known in the theory
of empirical processes (see, e.g., \cite{Dud} or \cite{Well}; we do not make countability assumption in
some of the examples below). 
Estimates $\hat f_n$ are supposed to be $\Sigma \times {\cal A}$-measurable. 
We denote by $P_n$ the empirical measure based on the sample $(X_1,\dots, X_n):$ 
$$
P_n := n^{-1}\sum_{j=1}^n \delta_{X_j},
$$
where $\delta_x$ is the probability measure concentrated at the point $x\in S.$
We also use the notation $\|\cdot \|_{\cal  F}$ for the sup-norm of functions 
from the class ${\cal F}$ into ${\Bbb R}:$ 
$$
\|Y\|_{\cal  F}:= \sup_{f\in {\cal  F}}|Y(f)|.
$$

Our approach is based on the following simple idea. 
Denote $B({r}):=\{f:P|f|\leq r\}$ and set $r_0^n=1.$
It's clear that for any consistent estimate $\hat f_n$ $P_n \hat f_n =0$ and, hence, 
$$P\hat f_{n}\leq 
P_n \hat f_n + \|P_{n}-P\|_{\cal  F}=
\|P_{n}-P\|_{\cal  F}
=\|P_n-P\|_{{\cal  F}\cap B({r_0^n})}=:r_1^n.$$
Therefore,
$\hat f_n\in {\cal F}\bigcap B({r_1^n}).$
It means that actually
$$
P\hat f_{n}\leq 
P_n \hat f_n + \|P_{n}-P\|_{{\cal  F}\cap B_{r_1^n}}=
\|P_{n}-P\|_{{\cal  F}\cap B_{r_1^n}}.
$$
We can repeat this recursive procedure infinitely many times.
Namely, if
$r_{k+1}^n:=\|P_{n}-P\|_{{\cal  F}\cap B({r_k^n})},$
then, by induction,  
$P\hat f_{n}\leq r_{k}^n$ for any natural $k.$
It is also clear that the sequence $\{r_k^n\}$
is nonincreasing Indeed, by a simple induction argument,
we have that $r_k^n\leq r_{k-1}^n$ implies that
$$
r_{k+1}^n = 
\|P_n-P\|_{{\cal  F}\cap B(r_k^n)}
\leq 
\|P_n-P\|_{{\cal  F}\cap B(r_{k-1}^n)}=r_k^n.
$$
Thus, the following proposition holds.

\begin{proposition} 
The sequence $\{r_k^n\}_{k\geq 1}$ is nonincreasing 
and for any consistent estimate $\hat f_n$ $P\hat f_n\leq \inf_{k\geq 0}r_k^n.$
\end{proposition}

The sequence $\{r_k^n\}_{k\geq 1}$ depends not only on the data; it also depends 
explicitly on the unknown distribution $P,$ so it can not be used for the purposes
of bounding the risk. However, there is a simple bootstrap type
approach that allows one to get around this difficulty.

The Rademacher process indexed by the function class ${\cal F}$ is defined as
$$R_{n}=\frac{1}{n}\sum_{i=1}^{n}{\eps}_{i}{\delta}_{X_{i}},$$
where $\{\eps_i\}$ is a Rademacher sequence (an i.i.d. sequence of random variables taking the values 
$+1$ and $-1$ with probability $1/2$ each) independent of $\{X_i\}.$
It has been used for a long time 
to obtain the bounds on the sup-norm of the empirical process indexed by
functions (in the so called symmetrization inequalities, see \cite{Well}).
Recently, Koltchinskii \cite{Kol1} (see also \cite{KAADP}) suggested to use $\|R_n\|_{\cal  F}$
as data-based measure of the accuracy of empirical approximation $\|P_n-P\|_{\cal  F}$
in learning problems and developed a version of structural risk minimization
in which the norms of Rademacher process play the role of data-dependent penalties.
Lozano \cite{Loz} compared this method of penalization with the method based
on VC-dimensions and the cross-validation method and found out that in the so called 
problem of the "intervals model selection" the Rademacher penalization performs 
better than other methods. Hush and Scovel (1999) used Rademacher norms to obtain
posterior performance bounds for machine learning. 
However, the "global" norm of Rademacher process does not 
allow one to recover the rate of convergence of the risk to $0$ in the 
case when $f_0\in {\cal F}$ (the so called zero error case). To address this problem,
we define below a sequence of localized norms of Rademacher process that majorizes 
the sequence $\{r_k^n\}$ defined above. 

Given $\eps>0,$ let $\bar \varphi$ be a (random) function defined by
$$\bar \varphi (r):=\bar K_1 \|R_{n}\|_{{\cal  F}\cap B_{2r}^e}+\bar K_2\sqrt{r\eps} +\bar K_3 \eps,$$
where $B_{r}^{e}=\{f\in {\cal F}:P_{n}f\leq r\}$ and $\bar K_1, \bar K_2, \bar K_3>0$ are numerical
constants. 

We introduce the following data-dependent sequence 
$$\{\bar r_k^n\}_{k\geq 0}=\{\bar r_k^n(X_1,\dots ,X_n;\eps_1,\dots ,\eps_n)\}_{k\geq 0},$$ 
\begin{equation}
{\bar r}_{0}^{n}=1,\,\,
{\bar r}_{k+1}^{n}=\bar \varphi({\bar r}_k^{n})\wedge 1, 
\,\,\,\,\,\,\,\,
k=0,1,2,\ldots \label{Seq}
\end{equation}
Since the function $\bar \varphi$ is nondecreasing, a simple induction shows that the sequence $\{\bar r_k^n\}$
is nonincreasing.

\begin{theorem} 
There is a choice of numerical constants $\bar K_1, \bar K_2, \bar K_3>0$ such that 
for all $P\in {\cal P}(S),$ for all $N\geq 1$ and for any consistent estimate $\hat f_n$
$$
{\Bbb P}\bigl\{P\hat f_n\geq \bar r_N^n  \bigr\}\leq 2 N e^{-\frac{n\eps}{2}}.
$$ 
\end{theorem}

Thus, if one chooses $N\geq 1$ and, for a given $\delta >0,$ 
$\eps >(\log {{2N\delta}})/{n},$
then one can define $\hat \beta_n ({\cal F};\delta):=\bar r_N^n$ to get the bound (1.1).
The question to be answered is how large should be the number of iterations $N$ to achieve 
a reasonably good upper bound on the risk in such a way (if it is possible at all). Surprisingly,
under rather general conditions the upper bound becomes sharp after very few 
iterations (roughly, the number of iterations $N$ is of the order $\log_2 \log_2 (\frac{1}{\eps})$). 

In what follows, given a (pseudo)metric space $(M;d),$ we denote $N_d(M;\eps)$ the minimal number of balls of
radius $\eps ,$ covering $M,$ and $H_d(M;\eps):=\log {N_d(M;\eps)}.$ Also, for a probability measure $Q$ on $(S,{\cal A}),$ 
$d_{Q,2}$ denotes the metric of the space $L_2(S;dQ).$

Given a class of functions ${\cal F},$
assume that
$$
{\Bbb E}_{\eps}\|
n^{-1/2}\sum_{i=1}^{n}{\eps}_{i}{\delta}_{X_{i}}
\|_{B^e(r)\cap {\cal  F}}
\leq \hat \psi_n(\sqrt{r})
$$
for some concave nondecreasing (random) function $\hat \psi_n.$
Usually the role of $\hat \psi_n$ will be played by
the random entropy integral
$$
\hat \psi_n(r)=
K\int\limits_0^{r}H_{d_{P_n,2}}^{1/2}({\cal F},u)du
$$
or by some further upper bound on the random entropy integral.
Let us denote by $\hat \delta_n := \hat \delta_n (X_1,\dots ,X_n)$ the solution of the equation
$$
\hat \delta_n=n^{-1/2}\hat \psi_n\bigl(\sqrt{\hat \delta_n}
\bigr).
$$
The following theorem gives the upper bound on the quantity $\bar r_{N}^n.$

\begin{theorem}
{\it
If the number of iterations is equal to
$
N=[\log_2\log_2\eps^{-1}]+1, 
$ 
then for some numerical constant
$c>0$ and for all $P\in {\cal P}(S)$  
$$
{\Bbb P}\left(
{\bar r}_N^n\geq c(\hat \delta_n\vee\eps)
\right)\leq 
([\log_2\log_2\eps^{-1}]+1)
e^{-\frac{n\eps}{2}}.
$$
}
\end{theorem}

{\bf Example 1. Learning a concept from a VC-class}. 
Consider the case of the concept learning, when ${\cal F}:=\{I_C: C\in {\cal C}\}.$
Given a sample $(X_1,\dots ,X_n)$ with unknown common distribution $P\in {\cal P}(S),$
we observe the labels $\{Y_j:=I_{C_0}(X_j):1\leq j\leq n\}$ for an unkown target concept
$C_0\in {\cal C}.$ An estimate $\hat C_n = \hat C_n ((X_1,Y_1),\dots ,(X_n,Y_n))$
of the target concept $C_0$ is called consistent iff $I_{\hat C_n}(X_j)=Y_j$
for all $j=1,\dots ,n.$ 
Let
$$
\Delta^{\cal  C}(X_1,\dots ,X_n):={\rm card}\bigl(\bigl\{C\cap \{X_1,\dots ,X_n\}: C\in {\cal C}\bigr\}\bigr).
$$
Then 
$$
\hat \psi_n (r) := K (\log {\Delta^{\cal  C}(X_1,\dots ,X_n)})^{1/2} r
$$
is an upper bound on the random entropy integral, 
which yields the value of $\hat \delta_n$
$$
\hat \delta_n = K^2 
{\frac{\log {\Delta^{\cal  C}}(X_1,\dots, X_n)}{n}}.
$$
Thus, with the same choice of $N$ we get for some numerical constant $c>0$ the bound 
$$
{\Bbb P}\left(
{\bar r}_N^n\geq c\bigl({\frac{\log {\Delta^{\cal  C}(X_1,\dots, X_n)}}{n}}\vee\eps\bigr)
\right)\leq 
([\log_2\log_2\eps^{-1}]+1)
e^{-\frac{n\eps}{2}}.
$$
Theorem 2 implies at the same time that for any consistent estimate $\hat C_n$ 
we have $P(\hat C_n\bigtriangleup C_0)\leq \bar r_N^n$ with probability at least
$1-2Ne^{-n\eps/2}.$  
This shows that for a VC-class of concepts ${\cal C}$ with VC-dimension $V({\cal C})$ 
the local Rademacher norm $\bar r_N^n$ (which, according to Theorem 2, is an upper bound on the risk of
consistent concepts $\hat C_n$) 
is bounded from above by the quantity $O(V({\cal C})\log n/n).$
Up to a logarithmic factor, this is the optimal (in a minimax sense) convergence rate 
of the generalization error to $0$ (see, e.g., \cite{Dev}).

Next we consider the conditions in terms of entropy with bracketing
$H_{[\ ]}({\cal F},\eps):=\log N_{[\ ]}({\cal F},\eps ).$ 
Here $N_{[\ ]}({\cal F},\eps)$ denotes the minimal number
of "brackets" $[f^{-},f^{+}]:=\{f: f^{-}\leq f \leq f^{+}\}$
with $d_{P,2}(f^{-},f^{+})\leq \eps $
($f^{-}, f^{+}$ being two measurable functions from $S$ into $[0,1],$ such that $f^{-}\leq f^{+}$).  
Let 
$$
\psi_{[\ ]}(r)=\int_{0}^{r}\left(H_{[\ ]}({\cal F},u)+1\right)^{1/2}du.
$$
and let $\delta_{[n]}=\delta_{[n]}(P)$ be the solution of the equation
$$
\delta_{[n]} = n^{-1/2}\psi_{[\ ]}(\sqrt{\delta_{[n]}}).
$$
Again, we set for some $\eps >0$
$
N: =[\log_2\log_2 \eps^{-1}]+1.
$
Then the following theorem holds.

\begin{theorem}
There exists a constant $c>0$ such that for all $P\in {\cal P}(S)$
$$
{\Bbb P}\left({\bar r}_N^n\geq c(\delta_{[n]}(P) \vee \eps)
\right)\leq 
([\log_2\log_2 \eps^{-1}]+1)e^{-\frac{n\eps}{2}}.
$$
\end{theorem}

In particular, if $H_{[\ ]}({\cal F}; u)=O(u^{-\gamma}),$ where $\gamma <2,$ then
$\psi_{[\ ]}(r)\asymp r^{1-\gamma/2}$ and $\delta_{[n]}\asymp n^{-\frac{2}{2+\gamma}}.$

{\bf Example 2. Learning a concept from a $d$-dimensional cube}. Let $S=[0,1]^d.$ We consider a problem of
estimation of a set (a concept) 
$C_0\subset [0,1]^d,$
based on the observations $(X_j,Y_j),\ j=1,\dots ,n,$ where $X_j,\ j=1,\dots ,n$ are i.i.d. points
in $[0,1]^d$ with common distribution $P$ and $Y_j:=I_{C_0}(X_j),\ j=1,\dots ,n.$ Such a model
frequently occurs in the problems of edge estimation in image analysis (see Mammen and Tsybakov (1995)).
Assume that the distribution $P$ has a density $p$ such that for some $B>0$
$$
B^{-1} \leq p(x) \leq B,\ x\in [0,1]^d. 
$$
Let ${\cal C}$ be a class of Borel subsets in $[0,1]^d$ such that ${\cal C}\ni C_0.$
Let $\lambda$ be the Lebesgue measure on $[0,1]^d.$ 
Denote $N_{I}({\cal C};\eps)$ the minimal number of brackets
$[C^{-},C^{+}]:=\{C: C^{-}\subset C\subset C^{+}\}$ with $\lambda (C^{+}\setminus C^{-})\leq \eps $
($C^{-}, C^{+}$ being two measurable subsets in $[0,1]^d$ such that $C^{-}\subset C^{+}$). 
Let $H_{I}({\cal C};\eps):=\log N_{I}({\cal C};\eps).$ This version of entropy 
with bracketing is often called "entropy with inclusion". 
We define
$$
\psi_{I}(r)=\int_{0}^{r}\left(H_{I}({\cal C},u)+1\right)^{1/2}du,
$$
and let $\delta_{n}^{I}=\delta_{n}^{I}(P)$ be the solution of the equation
$$
\delta_{n}^I = n^{-1/2}\psi_{I}(\sqrt{\delta_{n}^I}).
$$
If we have 
$$
H_{I}({\cal C};u)=O(u^{-\gamma}),
$$
then Theorem 4 easily implies that 
with some constant $c>0$  
$$
{\Bbb P}\left({\bar r}_N^n\geq c(\delta_{n}^I \vee \eps)
\right)\leq 
([\log_2\log_2 \eps^{-1}]+1)e^{-\frac{n\eps}{2}},
$$
where $\delta_{n}^I\asymp n^{-\frac{1}{1+\gamma}}.$
By Theorem 2, for any consistent estimate $\hat C_n$ of the set $C_0$ 
(i.e. such that $I_{\hat C_n}(X_j)=Y_j,\ j=1,\dots ,n$),
the quantity $\bar r_N^n$ is an upper bound (up to a constant)
on $\lambda (\hat C_n\bigtriangleup C_0).$

In particular, if ${\cal C}$ is the class of sets with $\alpha$-smooth boundary
in $[0,1]^d,$ then well known bounds on the bracketing entropy due to Dudley
(see e.g. Dudley (1999)) imply that $\gamma = {\frac{d-1}{\alpha}}$ and
$
\delta_{n}^I = n^{-\frac{\alpha}{d-1+\alpha}}.
$ 
Similarly, if ${\cal C}$ is the class of closed convex subsets of $[0,1]^d,$
the rate becomes 
$
\delta_{n}^I = n^{-\frac{2}{d+1}}.
$ 
It was shown by Mammen and Tsybakov (1995) that both rates are optimal in a 
minimax sense. 

The examples above show that the local Rademacher penalties (defined only based
on the data and using neither prior information about the underlying distribution,
nor the specific properties of the function class) can recover the optimal
convergence rates of the estimates in function learning problems.


\section{Proofs of the main results}

The proofs of the results are based on a version of Talagrand's concentration inequalities
for empirical processes, see \cite{Tal1}, \cite{Tal2}.
The version of the inequalities we are using, with explicit numerical values of the 
constants involved (that determine the values of the constants in our procedures, 
such as $\bar K_1, \bar K_2, \bar K_3$ above)
are due to Massart (1999). These inequalities are also very convenient 
for applications since the quantity $\sigma^2$ (the sup-norm of the variances, see below) 
they involve is very easy to bound. 
It should be also mentioned that
the idea to use Talagrand's concentration inequalities to bound
the risk in nonparametric estimation and, especially, in model
selection problems goes back to Birg\' e and Massart (see \cite{Bir1},
\cite{Bar1} and references therein).  

We formulate now Massart's inequality in a form convenient for our purposes.

\begin{theorem} 
Let 
${\cal F}$
be some countable family of real valued measurable functions,
such that 
$\|f\|_{\infty}\leq b<\infty$
for every
$f\in {\cal F}.$
Let $Z$ denote either
$\|P_{n}-P\|_{\cal  F}$
or
$\|R_{n}\|_{\cal  F}.$
Let
$\sigma^2=n\sup {\rm Var}(f(X_{1})).$
Then for any positive real number $x$
and $0<\gamma<1$
\begin{equation} 
{\Bbb P}(Z\geq (1+\gamma){\Bbb E}Z +
[ \sigma\sqrt{2kx}+k(\gamma)bx]/n)
\leq e^{-x},\label{Massart1}
\end{equation}
where $k$ and $k(\gamma)$ can be taken equal to
$k=4$ and $k(\gamma)=3.5+32\gamma^{-1}.$
Moreover, one also has
\begin{equation}
{\Bbb P}(Z\leq (1-\gamma){\Bbb E}Z - 
[\sigma\sqrt{2k'x}-k'(\gamma)bx]/n)
\leq e^{-x},\label{Massart2}
\end{equation}
where $k'=5.4$ and $k'(\gamma)=3.5+43.2\gamma^{-1}.$
\end{theorem}

{\bf Proof of Theorem 2}. 
Let for any fixed real positive number $r$
$$
\varphi_1(r)=
\|P_{n}-P\|_{{\cal  F}\cap B({r})}
$$
$$
\varphi_2(r)=
(1+\gamma){\Bbb E}
\|P_{n}-P\|_{{\cal  F}\cap B({r})}+
2\sqrt{r\eps}
+(1.75+16{\gamma}^{-1})\eps .
$$
\begin{eqnarray*}
\varphi_3(r)=
\frac{2(1+\gamma)}{1-{\gamma}'}
\biggl[
\|R_{n}\|_{{\cal  F}\cap B({r})}
&+&
\sqrt{5.4r\eps}+(1.75+21.6{{\gamma}'}^{-1})\eps
\biggr]
\\
& + &
2\sqrt{r\eps}
+(1.75+16{\gamma}^{-1})\eps .
\end{eqnarray*}
Then, for any $r>0$
\begin{equation}
{\Bbb P}\Bigl(\varphi_1(r)\leq \varphi_2(r) \leq \varphi_3(r)\Bigr)
\geq 1-2e^{-\frac{n\eps}{2}}.\label{lemma1}
\end{equation}

Indeed, in order to apply inequalities (\ref{Massart1}) 
and (\ref{Massart2}), we notice that for every
$f\in {\cal F}\bigcap B({r})$ the sup-norm 
$\|f\|_{\infty}\leq b=1$ and
$$
\sigma^2=\sup_{{\cal  F}\cap B_{r}}n{\rm Var}(f(X))
\leq \sup_{{\cal  F}\cap B({r})} n Pf^2
\leq \sup_{{\cal  F}\cap B({r})} n Pf
\leq nr.
$$ 
Moreover, if we set 
$x=n\eps/2$, 
then (\ref{Massart1}) implies
\begin{eqnarray*}
&&
{\Bbb P}\Bigl(\|P_{n}-P\|_{{\cal  F}\cap B({r})}
\geq  
(1+\gamma){\Bbb E} \|P_{n}-P\|_{{\cal  F}\cap B({r})}
+2\sqrt{r\eps}
\\
&&
+(1.75+16\gamma^{-1})\eps\Bigr)
\leq e^{-\frac{n\eps}{2}},
\end{eqnarray*}
and (\ref{Massart2}) implies
\begin{eqnarray*}
&&
{\Bbb P}\Bigl({\Bbb E}\|R_{n}\|_{{\cal  F}\cap B({r})}
\geq 
(1-\gamma')^{-1} [\|R_{n}\|_{{\cal  F}\cap B({r})}
+\sqrt{5.4r\eps}
\\
&&
+(1.75+21.6\gamma'{}^{-1})\eps]\Bigr)
\leq e^{-\frac{n\eps}{2}}.
\end{eqnarray*}
Taking into account the symmetrization inequality
$$
{\Bbb E} \|P_{n}-P\|_{{\cal  F}\cap B({r})}
\leq
2{\Bbb E}\|R_{n}\|_{{\cal  F}\cap B({r})},
$$ 
we get (2.3).

We set 
$$
\bar K_1 := \frac{2(1+\gamma)}{1-{\gamma}'},\ \bar K_2 := \frac{2\sqrt{5.4}(1+\gamma)}{1-{\gamma}'}
+ 2,
$$
$$
\bar K_3 := \frac{2(1+\gamma)}{1-{\gamma}'}(1.75+21.6{{\gamma}'}^{-1})+(1.75+16{\gamma}^{-1}). 
$$

Let us introduce the following sequence:
$ 
{\hat r}_0^n :=1
$
and
$
{\hat r}_{k+1}^n =\varphi_2({\hat r}_k^n)\wedge 1
$ 
for $k=0,1,2,\ldots .$
Since $\varphi_2$ is nondecreasing, it's easy to prove by induction
that the sequence $\{\hat r_k^n\}$ is nonincreasing. 

We will also prove by induction that for all $k\geq 0$ 
\begin{equation}
{\Bbb P}\Bigl\{r_i^n\leq \hat r_i^n \leq \bar r_i^n,\ i\leq k\Bigr\}\geq 
1-2ke^{-\frac{n\eps}{2}}.
\end{equation}
For $k=0$ (2.4) is trivial since $r_0^n=\hat r_0^n = \bar r_0^n =1.$
We proceed by the induction argument.
Let us introduce the events
$$
{\cal A}_k=\{r_i^n \leq {\hat r}_i^n \leq {\bar r}_i^n,\,\,i\leq k\}
\,\,\,\mbox{ and }
{\cal B}_k=\{\varphi_1({\hat r}_k^n) \leq \varphi_2({\hat r}_k^n) 
\leq \varphi_3({\hat r}_k^n)\}.
$$
To make the induction step, let us 
assume that we have already proven that
$$
{\Bbb P}\left({\cal A}_k\right)
\geq 1-2ke^{-\frac{n\eps}{2}}.
$$
Then (2.3) implies 
$$
{\Bbb P}\left({\cal B}_k\right)
\geq 1-2e^{-\frac{n\eps}{2}}.
$$
On the event ${\cal A}_k\bigcap {\cal B}_k,$
$$
{\cal F}\cap B({\hat r}_k^n)\subseteq 
{\cal F}\cap B^e(2{\hat r}_k^n),
$$
since for $f\in {\cal F}\bigcap B({{\hat r}_k^n})$
\begin{eqnarray*}
P_{n}f
&\leq& 
Pf + \|P_n-P\|_{{\cal  F}\cap B({\hat r}_k^n)}\leq 
{\hat r}_k^n
+\|P_n-P\|_{{\cal  F}\cap B({\hat r}_k^n)}
\\
&=&
{\hat r}_k^n+\varphi_1({\hat r}_k^n)
\leq
{\hat r}_k^n+\varphi_2({\hat r}_k^n)
={\hat r}_k^n+{\hat r}_{k+1}^n
\leq 2{\hat r}_k^n,
\end{eqnarray*}
which implies that the inequalities
$\varphi_3({\hat r}_k^n)\leq \bar \varphi({\hat r}_k^n) 
\leq \bar \varphi({\bar r}_k^n)={\bar r}_{k+1}^n$
hold.
Therefore, on the event ${\cal A}_k\bigcap {\cal B}_k,$
$$
r_{k+1}^n=\varphi_1(r_{k}^n)\leq \varphi_1({\hat r}_{k}^n)
\leq \varphi_2({\hat r}_{k}^n)= {\hat r}_{k+1}^n\leq
\varphi_3({\hat r}_{k}^n)\leq {\bar r}_{k+1}^n.
$$
So, ${\cal A}_k\bigcap {\cal B}_k\subseteq {\cal A}_{k+1},$
that completes the proof of the induction step
$$
{\Bbb P}\left({\cal A}_{k+1}\right)
\geq 1-2(k+1)e^{-\frac{n\eps}{2}}.
$$
It follows that 
$$
{\Bbb P}(r_N^n> {\bar r}_N^n)\leq 
2N e^{-\frac{n\eps}{2}},
$$
and since, by Proposition 1, $P\hat f_n \leq r_N^n,$ we conclude that 
$$
{\Bbb P}\{P\hat f_n > \bar r_N^n\}\leq 
2N e^{-\frac{n\eps}{2}}. 
$$

{\bf Proof of Theorem 3}. Let $(\Omega_{\eps}, {\Sigma}_{\eps}, {\Bbb P}_{\eps})$
denote the probability space on which the Rademacher sequence
$\eps_1,\ldots,\eps_n, \dots$ is defined, ${\Bbb E}_{\eps}$
being the expectation with respect to ${\Bbb P}_{\eps}.$  
We introduce the function
\begin{eqnarray}
&&
\varphi_4(r) 
 = 
\frac{2(1+\gamma)}{1-{\gamma}'}
\biggl[
(1+{{\gamma}''}^{-1}){\Bbb E}_{\eps}\|R_{n}\|_{{\cal  F}\cap B^e(2r)}
+2\sqrt{r\eps}
\nonumber
\\
&&
+
(1.75+16{{\gamma}''}^{-1})\eps
+\sqrt{5.4r\eps}+(1.75+21.6{{\gamma}'}^{-1})\eps
\biggr]
\nonumber
\\
&&
+
2\sqrt{r\eps}
+(1.75+16{\gamma}^{-1})\eps ,\label{phi4}
\end{eqnarray}
where
${\gamma}''>0.$
The inequalities (\ref{Massart1}) and (\ref{Massart2})
also hold for the conditional probability
${\Bbb P}_{\eps}$
and the process
$Z=R_n$
with fixed 
$X_1,\ldots, X_n.$
Therefore, for any
$r>0$

$$ 
{\Bbb P}_{\eps}(\bar \varphi(r)\leq \varphi_4(r))
\geq 1-e^{-\frac{n\eps}{2}}.
$$

Define a sequence
$$
{\check r}_0^n=\varphi_4(1),\,\,\,
{\check r}_{k+1}^n=\varphi_4({\check r}_{k}^n)\wedge 1 
,\,\,\,\,
k=0,1,2,\ldots
$$
By the induction argument, similar to the one we used
in the proof of theorem 2, we get
$$
{\Bbb P}_{\eps}\biggl(
\bigcap_{i=1}^{N}
\{{\bar r}_i^n\leq {\check r}_i^n\}
\biggr)
\geq 1-Ne^{\frac{n\eps}{2}}.
$$  
If we prove that
${\check r}_{k}^n\leq a_k$
for a sequence $a_k,$
independent of
$\eps_1,\ldots,\eps_n,$
then the unconditional probability
$$
{\Bbb P}\biggl(
\bigcap_{i=1}^{N}
\{{\bar r}_i^n\leq a_i\}
\biggr)
\geq 1-Ne^{\frac{n\eps}{2}}.
$$

By the assumption we have
\begin{equation}
{\Bbb E}_{\eps}\|
n^{-1}\sum_{i=1}^{n}{\eps}_{i}{\delta}_{X_{i}}
\|_{B^e(r)\cap {\cal  F}}
\leq
\hat \psi_n(\sqrt{r}).\label{Ebound}
\end{equation}

Hence, we can choose $c\geq 1$,
depending on the parameters
$\gamma,\gamma',\gamma''$
in the definition (\ref{phi4}) of the function
$\varphi_4,$
in such a way that
$$
{\check r}_{k+1}^n=\varphi_4 ({\check r}_k^n) \leq 
c\left(
\eps+({\check r}_k^n\eps)^{1/2}+
n^{-1/2}\hat \psi_n\left(\sqrt{{\check r}_k^n}\right)
\right).
$$
The above inequality implies by induction that the sequence
$$
r_0=1,\,\,
r_{k+1}=
c\left(
\eps+(r_k\eps)^{1/2}+
n^{-1/2}\hat \psi_n\left(\sqrt{r_k}\right)
\right)\wedge 1,
$$
majorizes the sequence ${\check r}_k^n.$

It's clear that in the case when
$r_1<1$
the sequence 
$r_k$ is decreasing and it converges to the solution 
$\delta$
of the equation
$$
\delta=
c\left(
\eps+(\delta\eps)^{1/2}+
n^{-1/2}\psi\left(\sqrt{\delta}\right)
\right).
$$
Let us study the behaviour of the difference
$
d_k:=r_k-\delta.
$
Since the function $\hat \psi_n$ is concave, we have 
$$
\hat \psi_n^{\prime}(\sqrt{\delta})\leq
\hat \psi_n(\sqrt{\delta})/\sqrt{\delta}.
$$
The definition of $\delta\,\,$
implies that
$$
c\left(
n^{-1/2}\hat \psi_n(\sqrt{\delta})+\sqrt{\delta\eps}
\right)\leq \delta.
$$ 
Therefore
\begin{eqnarray*}
& &
d_{k+1}
=
r_{k+1}-\delta=
c\left(
n^{-1/2}\hat \psi_n(\sqrt{r_k})-n^{-1/2}\hat \psi_n(\sqrt{\delta})+
\sqrt{r_k\eps}-\sqrt{\delta\eps}
\right)
\\
& &
\leq
c\left(
n^{-1/2}\hat \psi_n^{\prime}(\delta)+\sqrt{\eps}
\right)\sqrt{r_k-\delta}
\leq
c\left(
n^{-1/2}\hat \psi_n(\sqrt{\delta})+\sqrt{\delta\eps}
\right)/\sqrt{\delta}\sqrt{d_k}
\\
& &
\leq
\sqrt{\delta d_k}.
\end{eqnarray*}
We have proven that the sequence 
$d_k$ satisfies the following inequality 
$$
d_{k+1}\leq \sqrt{\delta d_k},\ k\geq 0.
$$
Now it's easy to show by induction that
$$
d_N\leq \delta^{2^{-1}+\ldots+2^{-N}}
= \delta^{1-2^{-N}}.
$$
Going back to the sequence $r_k,$ we get that
$$
r_N=\delta+d_N\leq
\delta\left(1+\delta^{-2^{-N}}\right).
$$
Since the definition of $\delta$ implies that
$\delta^{-1}<\eps^{-1},$ then the choice of 
$$
N=\left[
\log_2\log_2 \eps^{-1}
\right]+1
$$
guarantees that 
$\delta^{-2^{-N}}\leq 2$ and, hence, 
$
r_N\leq (1+2)\delta=3\delta.
$
What remains to do in order to finish the proof of the theorem, 
is to bound $\delta$ by the maximum of $\eps$
and the solution $\hat \delta_n$ of the equation
$
\hat \delta_n=n^{-1/2}\hat \psi_n(\sqrt{\hat \delta_n}).
$
Actually, we will prove that $\delta$ is bounded dy
$
\delta'':=(3c)^2 \delta',
$
where
$
\delta'=
\left(
\hat \delta_n\vee \eps.
\right)
$
First of all let us notice that the fact
that $\hat \psi_n$ is concave and $\hat \psi_n(0)=0$
implies that for $c\geq 1$
$\hat \psi_n(cx)\leq c\hat \psi_n(x).$
Also note that, since $\delta' \geq \hat \delta_n,$ the concavity of $\hat \psi_n$ 
and the definition of $\hat \delta_n$ imply  
$$
n^{-1/2}\hat \psi_n\left(\sqrt{\delta'}\right)\leq
\frac{n^{-1/2}\hat \psi_n(\sqrt{\hat \delta_n})}{\sqrt{\hat \delta_n}}\sqrt{\delta'}
=\sqrt{\hat \delta_n} \sqrt{\delta'}\leq 
\delta'.
$$
Combining these properties, we get 
\begin{eqnarray*}
&&
c\left(
\eps+(9c^2\delta'\eps)^{1/2}+
n^{-1/2}\hat \psi_n\left(3c\sqrt{\delta'}\right)
\right)
\\
&&
\leq
c\left(2
\sqrt{(3c)^2}\delta'+
\delta'
\right)\leq
9c^2\delta'=\delta''.
\end{eqnarray*}
With necessity it means that 
$
\delta\leq \delta''=9c^2(\hat \delta_n\vee \eps).
$
And, hence, 
$
{\bar r}_N^n
\leq \delta''\leq 27c^2(\hat \delta_n\vee \eps).
$

The theorem is proven.

{\bf Proof of Theorem 4.}
In order to bound
${\bar r}_k,$ we first construct the bound on 
$\|R_n\|_{{\cal  F}\cap B^{e}(2{\bar r}_k)}$
in terms of
${\Bbb E}\|P_n-P\|_{{\cal  F}\cap B({\check r}_k)}$
for properly defined sequence 
${\check r}_k.$
Afterwards, the expectation can be majorized by
the bracketing entropy integral. We will show that the
sequence ${\check r}_k$
can be chosen as follows
$$
{\check r}_0=1,\,\,\,
{\check r}_{k+1}=\left(
\tilde c_1{\Bbb E}\|P_n-P\|_{{\cal  F}\cap B(3{\check r}_k)}+
\tilde c_2\sqrt{\eps {\check r}_k }+\tilde c_3
\right)\wedge 1,
$$
for some large enough constants 
$\tilde c_1, \tilde c_2, \tilde c_3>0.$ 
One can argue similarly to the proof of Theorem 3 to show that
the following bound holds: 
\begin{equation}
{\Bbb P}\Bigl(
\bigcap_{k\leq i}\{
{\bar r}_k\leq {\check r}_k
\}
\Bigr)\geq 1-2ie^{-\frac{n\eps}{2}}.
\end{equation}
We will prove even a stronger assertion 
that for the event
$$
{\cal A}_i=
\bigcap_{k\leq i}\Bigl(\{
{\bar r}_k\leq {\check r}_k
\}\cap\{{\cal F}\cap B^{e}(2{\bar r}_k)\subseteq {\cal F}\cap B(3{\check r}_k)\}\Bigr)
$$ 
we have
\begin{equation}
{\Bbb P}({\cal A}_i)\geq 1-2ie^{-\frac{n\eps}{2}}.
\end{equation}

Let us choose the constants 
$c_1^{\prime}, c_2^{\prime}, c_3^{\prime }>0$ and $\tilde c_1, \tilde c_2, \tilde
c_3>0$ in such a way 
that for the functions
$$
\varphi_{5}(r)=\left(
c_1^{\prime}\|P_n-P\|_{{\cal  F}\cap B(r)}+
c_2^{\prime}\sqrt{\eps r} + c_3^{\prime}\eps
\right)
$$
and
$$
\varphi_{6}(r)=\left(
\tilde c_1{\Bbb E}\|P_n-P\|_{{\cal  F}\cap B(r)}+
\tilde c_2\sqrt{\eps r}+\tilde c_3 \eps 
\right),
$$
the inequalities of Massart (see Theorem 5) would imply that for any fixed
$r>0$
$$
\varphi_3(r)\leq \varphi_5(r)\leq\varphi_6(r)
$$
with probability at least
$1-2e^{-\frac{n\eps}{2}}$
(the function $\varphi_3$ was defined in the proof of Theorem 2). 
Clearly, we have 
$
\check r_{k+1} = \varphi_6 (\check r_k)\wedge 1.
$

First observe that (2.8) holds for $i=0$ (since
$\bar r_0=\check r_0=1$).  
Define
$$
{\cal B}_i:=\{\varphi_3(3\check r_i)\leq \varphi_5(3\check r_i)\leq \varphi_6(\check r_i)\}.
$$
Then 
$$
{\Bbb P}({\cal B}_i) \geq 
1-2e^{-\frac{n\eps}{2}}.
$$
To make an induction step, we first of all
notice that on the event 
${\cal A}_i\cap {\cal B}_i,$ we have
$$
{\bar r}_{i+1}=\bar \varphi (\bar r_i)\wedge 1 \leq 
\varphi_3(3{\check r}_i)\wedge 1\leq
\varphi_5(3{\check r}_i)\wedge 1\leq
\varphi_6(3{\check r}_i)\wedge 1= {\check r}_{i+1}.
$$
Also, on the event ${\cal A}_i\cap {\cal B}_i,$ we have
$ 
{\cal F}\cap B^{e}(2{\bar r}_{i+1})\subseteq {\cal F}\cap B(3{\check r}_{i+1}).
$
Indeed, if 
$f\in {\cal F}\cap B^{e}(2{\bar r}_{i+1}),$ then
\begin{eqnarray*}
Pf &\leq & 2{\bar r}_{i+1}+
\|P_n -P\|_{{\cal  F}\cap B^{e}(2{\bar r}_{i+1})}\leq
2{\bar r}_{i+1}+
\|P_n -P\|_{{\cal  F}\cap B^{e}(2{\bar r}_{i})}
\\
&\leq&
2{\bar r}_{i+1}+
\|P_n -P\|_{{\cal  F}\cap B(3{\check r}_{i})}\leq
2{\bar r}_{i+1}+
\varphi_5(3{\check r}_i)\wedge 1
\\
&\leq&
2\bar r_{i+1}+\varphi_6(3{\check r}_i)\wedge 1
= 2\bar r_{i+1}+\check r_{i+1} 
\leq 3{\check r}_{i+1}
\end{eqnarray*}
(to show that $\|P_n-P\|_{{\cal  F}\cap B(3\check r_i)}\leq \varphi_5(3\check r_i)\wedge 1$ 
we used the fact that the costant $c_1^{\prime}$ in the definition of ${\varphi}_5$
is larger than $1$). 
Thus, ${\cal A}_i\cap {\cal B}_i\subset {\cal A}_{i+1}$ and
$$
{\Bbb P}({\cal A}_{i+1})\geq 
1-2(i+1) e^{-\frac{n\eps}{2}}.
$$
The proof of the induction step and of the bounds (2.8) and (2.7) is complete.

To finish the proof of the theorem
one has to bound
${\Bbb E}\|P_n-P\|_{{\cal  F}\cap B(r)}.$
Since for all $g\in {\cal F}\bigcap B(r)$ 
we have 
$\|g\|_{P,2}\leq (Pg)^{1/2}\leq \sqrt{r}$
and $|g|\leq 1$ then by Theorem 2.14.2
in \cite{Well} 
$$
{\Bbb E}\|P_n-P\|_{{\cal  F}\cap B(r)}\leq
c\left(
n^{-1/2}\psi_{[\ ]}\left(
\sqrt{r}\right)
+I\{1>\sqrt{n}a(\sqrt{r})\}
\right),
$$
where
$$
a(\sqrt{r})=\sqrt{r}/
\sqrt{1+H_{[\ ]}({\cal F},\sqrt{r})}.
$$
We can assume that 
${\check r}_N\geq \delta_{[n]},$
otherwise, bound (2.7) immediately implies the 
assertion of the theorem. Therefore, 
${\check r}_k\geq \delta_{[n]}$ for all $k\leq N,$
which implies that 
$1\leq\sqrt{n}a(\sqrt{3{\check r}_k}).$
Indeed, using concavity of $\psi_{[\ ]}$ and the definition of $\delta_{[n]},$ we have
$$
\frac{\psi_{[\ ]}(\sqrt{3\check r_k})}{\sqrt{3\check r_k}}\leq 
\frac{\psi_{[\ ]}(\sqrt{\delta_{[n]}})}{\sqrt{\delta_{[n]}}}=
\sqrt{n}\sqrt{\delta_{[n]}}\leq \sqrt{n}\sqrt{3\check r_k}, 
$$
which implies 
$$
3{\check r}_k\geq n^{-1/2}\psi_{[\ ]}\left(\sqrt{3{\check r}_k}\right)
\geq n^{-1/2}\left(3{\check r}_k\right)^{1/2}
\left(1+H_{[\ ]}\left({\cal F}, \sqrt{3{\check r}_k}\right)
\right)^{1/2}.
$$
Hence,
$1\leq\sqrt{n}a(\sqrt{3{\check r}_k})$
and
$$
{\Bbb E}\|P_n-P\|_{{\cal  F}\cap B(3{\check r}_k)}\leq
cn^{-1/2}\psi_{[\ ]}\left(
\sqrt{3{\check r}_k}\right).
$$
Finally, with some constant $c>0$
$$
{\check r}_{k+1}\leq c\left(
n^{-1/2}\psi_{[\ ]}(\sqrt{3{\check r}_k})+\eps+
\sqrt{\eps {\check r}_k}
\right).
$$
The proof can be completed by the argument we used in Theorem 3.

\vskip 1mm

{\bf Acknowledgement}. 
The research of V. Koltchinskii is partially supported by NSA Grant
MDA904-99-1-0031. The research of D. Panchenko was 
partially supported by Boeing Computer Services Grant 3-48181.
The authors are thankful to Jon Wellner for pointing out the recent 
paper of Massart (1998).

\vskip 2mm

Department of Mathematics and Statistics

The University of New Mexico

Albuquerque NM 87131-1141

e-mail: vlad@math.unm.edu; panchenk@math.unm.edu

\end{document}